\documentclass[12pt,reqno]{article}

\usepackage[usenames]{color}
\usepackage{amssymb}
\usepackage{amsmath}
\usepackage{amsthm}
\usepackage{amsfonts}
\usepackage{amscd}

\usepackage[colorlinks=true,
linkcolor=webgreen,
filecolor=webbrown,
citecolor=webgreen]{hyperref}

\definecolor{webgreen}{rgb}{0,.5,0}
\definecolor{webbrown}{rgb}{.6,0,0}

\newcommand{\bell}{\textup{B}}
\usepackage{color}
\usepackage{fullpage}
\usepackage{float}
\usepackage{latexsym}
\usepackage{epsf}
\usepackage{breakurl}

\begin{document}

\theoremstyle{plain}
\newtheorem{theorem}{Theorem}
\newtheorem{corollary}[theorem]{Corollary}
\newtheorem{lemma}[theorem]{Lemma}
\newtheorem{proposition}[theorem]{Proposition}

\theoremstyle{definition}
\newtheorem{definition}[theorem]{Definition}
\newtheorem{example}[theorem]{Example}
\newtheorem{conjecture}[theorem]{Conjecture}

\theoremstyle{remark}
\newtheorem{remark}[theorem]{Remark}

\begin{center}
\vskip 1cm{\LARGE\bf 
A closed form for the generalized Bernoulli polynomials via Fa\`{a} di Bruno's formula
}
\vskip 1cm
\large
Sumit Kumar Jha\\
International Institute of Information Technology\\
Hyderabad-500 032, India\\
\href{mailto:kumarjha.sumit@research.iiit.ac.in}{\tt kumarjha.sumit@research.iiit.ac.in}\\
\end{center}

\vskip .2 in

\begin{abstract}
We derive a closed form for the generalized Bernoulli polynomial of order $a$ in terms of Bell polynomials and Stirling numbers of the second kind using the Fa\`{a} di Bruno's formula. 
\end{abstract}
\section{Main result}
\begin{definition}
The \emph{generalized Bernoulli polynomial} of order $a$ ($a\in \mathbb{N}$), denoted by $B_{n}^{a}(x)$, can be defined by the following exponential generating function
\begin{equation}
\left(\frac{t}{e^{t}-1}\right)^{a}e^{-xt}=\sum_{k=0}^{\infty}\frac{B^{a}_{k}(-x)}{k!}t^{k},
\end{equation}
where $|t|<2\pi$.\par 
Two closed forms for $B_{n}^{a}(x)$ are given in \cite{Srivastava} and \cite{Choi}. We prove the following formula. 
\begin{theorem}
We have
\begin{equation}
\label{main}
B_{n}^{a}(-x)=\sum_{k=0}^{n}(-1)^{k}\, k!\, B_{n,k}(\lambda_1,\cdots,\lambda_{n-k+1})
\end{equation}
where
$$
\lambda_{m}=\sum_{l=0}^{m}S(l+a,a)x^{m-l}\binom{m}{l}\binom{l+a}{a}^{-1} \qquad(m=1,2,\cdots,n-k+1),
$$
and $B_{n,k}(\lambda_1,\lambda_2,\dotsc,\lambda_{n-k+1})$ are the Bell polynomials \cite[p. 206]{ Comtet} defined by
\begin{equation}
\bell_{n,k}(\lambda_1,\lambda_2,\dotsc,\lambda_{n-k+1})=\sum_{\substack{1\le i\le n,\ell_i\in\mathbb{N}\\ \sum_{i=1}^ni\ell_i=n\\ \sum_{i=1}^n\ell_i=k}}\frac{n!}{\prod_{i=1}^{n-k+1}\ell_i!} \prod_{i=1}^{n-k+1}\Bigl(\frac{\lambda_i}{i!}\Bigr)^{\ell_i},
\end{equation}
and $S(n,k)$ are the Stirling numbers of the second kind for $n\ge k\ge1$ which can be generated by
\begin{equation}\label{2stirling-gen-funct-exp}
\frac{(e^x-1)^k}{k!}=\sum_{n=k}^\infty S(n,k)\frac{x^n}{n!}, \quad k\in\mathbb{N}.
\end{equation}
\end{theorem}
\end{definition}
\begin{proof}
Let $f(t)=1/t$ and $g(t)=\left(\frac{e^{t}-1}{t}\right)^{a}e^{xt}$. Using Fa\`{a} di Bruno's formula \cite[p. 134]{Comtet}
$$
{d^n \over dx^n} f(g(t)) = \sum_{k=0}^n f^{(k)}(g(t))\cdot B_{n,k}\left(g'(t),g''(t),\dots,g^{(n-k+1)}(t)\right),
$$
and the expansion
\begin{align*}
g(t)=\left(\frac{e^{t}-1}{t}\right)^{a}e^{xt}=\frac{a!}{t^{a}}\sum_{l=a}^{\infty}\frac{S(l,a)\, t^{l}}{l!}\sum_{k=0}^{\infty}\frac{x^{k}t^{k}}{k!}\\
=\frac{a!}{t^{a}}\sum_{l=0}^{\infty}\frac{S(l+a,a)\, t^{l+a}}{(l+a)!}\sum_{k=0}^{\infty}\frac{x^{k}t^{k}}{k!}=a!\, \sum_{l=0}^{\infty}\frac{S(l+a,a)\, t^{l}}{(l+a)!}\sum_{k=0}^{\infty}\frac{x^{k}t^{k}}{k!}\\
=a!\,\sum_{m=0}^{\infty}t^{m}\sum_{k+l=m}\frac{S(l+a,a)x^{k}}{(l+a)!\, k!}\\
=a!\, \sum_{m=0}^{\infty}t^{m} \sum_{l=0}^{m}\frac{S(l+a,a)\, x^{m-l}}{(l+a)!\, (m-l)!}\\
=a!\, \sum_{m=0}^{\infty}\frac{t^{m}}{m!}\sum_{l=0}^{m}\frac{S(l+a,a)\, x^{m-l}\, m!}{(l+a)!\, (m-l)!}\\
=\sum_{m=0}^{\infty}\frac{t^{m}}{m!}\sum_{l=0}^{m}\frac{S(l+a,a)\, x^{m-l}\, m!}{\frac{(l+a)!}{l!\, a!}\, (m-l)!\,l!}\\
=\sum_{m=0}^{\infty}\frac{t^{m}}{m!}\sum_{l=0}^{m}S(l+a,a)x^{m-l}\binom{m}{l}\binom{l+a}{a}^{-1}=\sum_{k=0}^{\infty}\frac{g^{(k)}(0)\, t^{k}}{k!}.
\end{align*}
we conclude our formula \eqref{main}. 
\end{proof}
\begin{remark}
Letting $x=0$ and $a=1$ in equation \eqref{main} gives us
$$
B_{n}=\sum_{k=0}^{n}(-1)^{k}\, k!\, B_{n,k}\left(\frac{1}{2},\frac{1}{3},\cdots,\frac{1}{n-k+1}\right)
$$
which was obtained in \cite{Qi}. Here $B_{n}$ denote the Bernoulli numbers.
\end{remark}
\begin{remark}
Similar method yields a closed form for the Bernoulli polynomials $B_{n}(x)$ obtained in \cite{Qi2}.
\subsection*{Evaluation of $B_{n,k}(\lambda_1,\cdots,\lambda_{n-k+1})$}
We have
\begin{align*}
{\displaystyle \sum _{n=k}^{\infty }B_{n,k}(\lambda_{1},\ldots ,\lambda_{n-k+1}){\frac {t^{n}}{n!}}= {\frac {1}{k!}}\left(\sum _{j=1}^{\infty }\lambda_{j}{\frac {t^{j}}{j!}}\right)^{k} }\\
=\frac{1}{k!}(g(t)-1)^{k}=\frac{1}{k!}\sum_{r=0}^{k}(-1)^{k-r}\binom{k}{r}g(t)^{r}\\
=\frac{1}{k!}\sum_{r=0}^{k}(-1)^{k-r}\binom{k}{r}\sum_{n=0}^{\infty}\frac{t^{n}}{n!}\sum_{l=0}^{n}S(l+ar,ar)(rx)^{n-l}\binom{n}{l}\binom{l+ar}{ar}^{-1}.
\end{align*}
Thus
\[
B_{n,k}(\lambda_1,\cdots, \lambda_{n-k+1})=\frac{1}{k!}\sum_{r=0}^{k}(-1)^{k-r}\binom{k}{r}\sum_{l=0}^{n}S(l+ar,ar)(rx)^{n-l}\binom{n}{l}\binom{l+ar}{ar}^{-1}.
\]
\end{remark}
We can use the above to conclude from equation \eqref{main} that
\begin{align*}
B_{n}^{a}(-x)=\sum_{k=0}^{n}\sum_{r=0}^{k}(-1)^{r}\binom{k}{r}\sum_{l=0}^{n}S(l+ar,ar)(rx)^{n-l}\binom{n}{l}\binom{l+ar}{ar}^{-1}\\
=\sum_{r=0}^{n}\sum_{k=r}^{n}(-1)^{r}\binom{k}{r}\sum_{l=0}^{n}S(l+ar,ar)(rx)^{n-l}\binom{n}{l}\binom{l+ar}{ar}^{-1}\\
=\sum_{r=0}^{n}\sum_{k=r}^{n}\binom{k}{r}\, (-1)^{r}\sum_{l=0}^{n}S(l+ar,ar)(rx)^{n-l}\binom{n}{l}\binom{l+ar}{ar}^{-1}\\
=\sum_{r=0}^{n}\binom{n+1}{r+1}\,(-1)^{r}\sum_{l=0}^{n}S(l+ar,ar)(rx)^{n-l}\binom{n}{l}\binom{l+ar}{ar}^{-1}.
\end{align*}

\end{document}